\documentclass[12pt]{article}

\usepackage{amsmath}
\usepackage{amsfonts}
\usepackage{bigints}

\begin{document}

\renewcommand{\P}{\par \rm Proof:\ }
\newcommand{\Pe}{$\hfill{\Box}$\break}
\newcommand{\Pes}{$\hfill{\Box}$}
\newcommand{\A}{\mbox{${{{\cal A}}}$}}


\author{Attila Losonczi}
\title{Means of infinite sets II}

\date{5 August 2018}

\newtheorem{thm}{\qquad Theorem}[section]
\newtheorem{prp}[thm]{\qquad Proposition}
\newtheorem{lem}[thm]{\qquad Lemma}
\newtheorem{cor}[thm]{\qquad Corollary}
\newtheorem{rem}[thm]{\qquad Remark}
\newtheorem{ex}[thm]{\qquad Example}
\newtheorem{prb}[thm]{\qquad Problem}
\newtheorem{df}[thm]{\qquad Definition}

\maketitle

\begin{abstract}

\noindent

We continue the study of how one can define means of infinite sets. We introduce many new properties, investigate their relations to each other and how they can typify a mean. We collect the properties in property groups e.g. for monotonicity and continuity because there is no single way to define such notions, instead there is a wide variety.

\noindent
\footnotetext{\noindent
AMS (2010) Subject Classifications: 40A05, 40G05, 28A10, 28A78 \\

}

\end{abstract}

\section{Introduction}
In this paper we are going to continue the investigations started in \cite{lamis}. For basic definitions, examples, ideas, intentions please consult \cite{lamis}.

In this paper our aim is to create many new properties that may typify a mean in many ways. We also create many property groups for internality, functionality, monotonicity and continuity concepts because it turns out that those notions cannot be grabbed in one single way.
In the first section we also enumerate some possible properties of the domain of a mean on infinite sets.

\medskip
We do not restrict ourself to arithmetic type means any more. Instead we are going to begin to build the theory of means that are defined on some subset of $P(\mathbb{R})$ and satisfies some very basic expected properties. Nevertheless in most of the cases we apply properties which only arithmetic type means own.

All sets considered in this paper are bounded.

\subsection{Basic notations}
For easier readability we copy the basic notations from \cite{lamis}.

Throughout this paper function $\A()$ will denote the arithmetic mean of any number of variables. Moreover if $(a_n)$ is an infinite sequence and $\lim_{n\to\infty}\A(a_1,\dots,a_n)$ exists then $\A((a_n))$ will denote its limit.

\begin{df}\label{davg}Let $\mu^s$ denote the s-dimensional Hausdorff measure ($0\leq s\leq 1$). If $0<\mu^s(H)<+\infty$ (i.e. $H$ is an $s$-set) and $H$ is $\mu^s$ measurable then $$Avg(H)=\frac{\int\limits_H x\ d\mu^s}{\mu^s(H)}.$$
If $0\leq s\leq 1$ then set $Avg^s=Avg|_{\{\text{measurable s-sets}\}}$. E.g. $Avg^1$ is $Avg$ on all Lebesgue measurable sets with positive measure.
\end{df}

\begin{df}For $K\subset\mathbb{R},\ y\in\mathbb{R}$ let us use the notation $$K^{-y}=K\cap(-\infty,y],K^{+y}=K\cap[y,+\infty).$$
\end{df}

Let $T_s$ denote the reflection to point $s\in\mathbb{R}$ that is $T_s(x)=2s-x\ (x\in\mathbb{R})$.

If $H\subset\mathbb{R},x\in\mathbb{R}$ then set $H+x=\{h+x:h\in H\}$. Similarly $\alpha H=\{\alpha h:h\in H\}\ (\alpha\in\mathbb{R})$.

We use the convention that these operations $+,\cdot$ have to be applied prior to the set theoretical operations, e.g. $H\cup K\cup L+x=H\cup K\cup (L+x)$.

$cl(H), H'$ will denote the closure and accumulation points of $H\subset\mathbb{R}$ respectively. Let $\varliminf H=\inf H',\ \varlimsup H=\sup H'$ for infinite bounded $H$.

Usually ${\cal{K}},{\cal{M}}$ will denote means, $Dom({\cal{K}})$ denotes the domain of ${\cal{K}}$.
\section{Properties revisited}
In this section we recall some properties already mentioned in \cite{lamis}, add many new ones to them and derive new relations among them. We also group the properties as there are some well defined property types.
\subsection{Possible properties of $Dom({\cal{K}})$}
We enumerate some possible properties of $Dom({\cal{K}})$. Only the first one is a preference from a mean however in most of the cases most of the mentioned properties will be fulfilled.
\begin{enumerate}
\item $Dom({\cal{K}})$ is closed under finite union and intersection. 

Although it seems very natural ${\cal{M}}^{iso}$ does not satisfies it, ${\cal{M}}^{iso}$ is not closed under finite union.
\item $H\in Dom({\cal{K}}), I$ is an interval then $H\cap I\in Dom({\cal{K}})$ if $H\cap I\ne\emptyset$.
\item $H^{-x}\in Dom({\cal{K}}),\ H^{+x}\in Dom({\cal{K}})$ then $H\in Dom({\cal{K}})$.

It is a consequence of 1.
\item $Dom({\cal{K}})$ is closed under translation, reflection and contraction/dilation. 

The arithmetic type means satisfy this one.
\item $H\in Dom({\cal{K}})$ then $cl(H)\in Dom({\cal{K}})$.

Usually means definied only on countable sets do not fulfill this one.
\item $H\in Dom({\cal{K}}), I$ is an interval then $I-H\in Dom({\cal{K}})$.

Means definied only on countable sets do not fulfill this one either.
\item $H\in Dom({\cal{K}})$ then $H'\in Dom({\cal{K}})$.
\item $Dom({\cal{K}})$ is closed under countable union/intersection.
\item If $H\in Dom({\cal{K}})$ then $f(H)\in Dom({\cal{K}})$ where $f:\mathbb{R}\to\mathbb{R}$ is any continuous function.
\end{enumerate}

\subsection{Properties related to internality}\label{ssexpprop}
\begin{itemize}
\item The very bacis property of a mean ${\cal{K}}$ is \textbf{internality} that is $$\inf H\leq {\cal{K}}(H)\leq\sup H.$$
It can happen that $\forall h\in H\ {\cal{K}}(H)<h$ (clearly when ${\cal{K}}(H)=\inf H$ and $\inf H\not\in H$) or in the opposite direction as well.

\item However almost always we require the stroger condition called \textbf{strong internality} $$\varliminf H\leq {\cal{K}}(H)\leq\varlimsup H.$$ This obviously implies for a set with only one accumulation point that its mean has be equal to the only accumulation point.

\item $\cal{K}$ has property \textbf{strict strong internality} if it is strongly internal and $\varliminf H<{\cal{K}}(H)<\varlimsup H$ whenever $H$ has at least 2 accumulation points.

\end{itemize}

\begin{prp}Clearly strict strong internal $\Rightarrow$ strong internal $\Rightarrow$ internal. \Pes
\end{prp}

\begin{prp}$Avg$ is strict strong internal.
\end{prp}
\P Let $H$ be a bounded $s$-set ($0\leq s\leq 1$). Let $a=\varliminf H,b=\varlimsup H$. Obviously there is $c\in (a,b)$ such that $\mu^s(H\cap[c,b])>0$. Let $H_1=H\cap[a,c),H_2=H\cap[c,b]$. Then $Avg(H)\geq \frac{\mu^s(H_1)a+\mu^s(H_2)c}{\mu^s(H)}=\frac{\mu^s(H_1)}{\mu^s(H)}a+\frac{\mu^s(H_2)}{\mu^s(H)}c>a$ because it is a weighted average and $\frac{\mu^s(H_2)}{\mu^s(H)}>0$.

The other inequality can be shown similarly.
\Pes

\begin{prp}${\cal{M}}^{lis}$ is strict strong internal. \Pes
\end{prp}

\begin{prp}$LAvg,{\cal{M}}^{iso},{\cal{M}}^{acc},{\cal{M}}^{eds}$ are not strict strong internal. \Pes
\end{prp}

\begin{ex}Example of a mean that is not strong internal:
$${\cal{K}}(H)=\frac{\inf\{\frac{x+y}{2}:x,y\in H\}+\sup\{\frac{x+y}{2}:x,y\in H\}}{2}.$$
\end{ex}

\subsection{Functional invariant properties}
\begin{df}If ${\cal{K}}$ is a mean then set $$F_{\cal{K}}=\{f:\mathbb{R}\to\mathbb{R}:H\in Dom({\cal{K}})\Rightarrow f(H)\in Dom({\cal{K}}),\ f({\cal{K}}(H))={\cal{K}}(f(H))\}.$$
\end{df}

Let us enumerate some functional invariant properties.

\begin{itemize}
\item The mean is \textbf{translation-invariant} if $x\in\mathbb{R}, H\in Dom({\cal{K}})$ then $H+x\in Dom({\cal{K}}),\ {\cal{K}}(H+x)={\cal{K}}(H)+x$. I.e. translations are in $F_{\cal{K}}$.

\item ${\cal{K}}$ is \textbf{point-symmetric} if $H\in Dom({\cal{K}})$ bounded and symmetric ($\exists s\in\mathbb{R}\ T_s(H)=H$) implies ${\cal{K}}(H)=s$.

\item ${\cal{K}}$ is \textbf{homogeneous} if $H\in Dom({\cal{K}}),\alpha\in\mathbb{R}$ then $\alpha H\in Dom({\cal{K}}),\ {\cal{K}}(\alpha H)=\alpha {\cal{K}}(H)$. I.e. contractions, dilations are in $F_{\cal{K}}$.

\item ${\cal{K}}$ is \textbf{reflection-invariant} if $H\in Dom({\cal{K}}),\ s\in\mathbb{R}$ then $T_s(H)\in Dom({\cal{K}}),\ {\cal{K}}(T_s(H))=T_s({\cal{K}}(H))$. I.e. reflections are in $F_{\cal{K}}$.

\end{itemize}

\begin{prp}If ${\cal{K}}$ reflection invariant then it is point-symmetric. \Pes
\end{prp}

\begin{prp}If ${\cal{K}}$ reflection invariant, $H\in Dom({\cal{K}}), {\cal{K}}(H)=s$ then ${\cal{K}}(T_s(H))=s$. \Pes
\end{prp}

\begin{prp}If $\cal{K}$ is point-symmetric, $s\in\mathbb{R},\ H\in Dom({\cal{K}})$ then ${\cal{K}}(H\cup T_s(H))=s$. \Pes
\end{prp}

\begin{prp}(a) If $f,g\in F_{\cal{K}}$ then $g\circ f \in F_{\cal{K}}$.

(b) If $f\in F_{\cal{K}},\ f$ has inverse, $H\in Dom({\cal{K}})$ implies $f^{-1}(H)\in Dom({\cal{K}})$ then $f^{-1} \in F_{\cal{K}}$.
\end{prp}
\P (a) If $H\in Dom({\cal{K}})$ then $f(H)\in Dom({\cal{K}})$ and $g({\cal{K}}(f(H)))={\cal{K}}(g(f(H)))={\cal{K}}(g\circ f(H))$ but $g({\cal{K}}(f(H)))=g(f({\cal{K}}(H)))=g\circ f({\cal{K}}(H))$.

(b) If $f^{-1}(H)\in Dom({\cal{K}})$ then $f({\cal{K}}(f^{-1}(H)))={\cal{K}}(f(f^{-1}(H)))={\cal{K}}(H)$ which gives that ${\cal{K}}(f^{-1}(H))=f^{-1}({\cal{K}}(H))$.
\Pes

\begin{ex}If we want to create a mean ${\cal{K}}$ that does not fulfill these properties then the easiest way is to construct a non-arithmetic type mean, e.g. let $f(x)=x^3,\ {\cal{K}}(H)=f^{-1}(Avg(f(H)))\ (H\in Dom(Avg))$. 
\end{ex}
\subsection{Monotonicity concepts}
\begin{itemize}
\item $\cal{K}$ is \textbf{monotone} if $\sup H_1\leq\inf H_2$ implies that ${\cal{K}}(H_1)\leq {\cal{K}}(H_1\cup H_2)\leq {\cal{K}}(H_2)$. 

${\cal{K}}$ is \textbf{strong monotone} if $\cal{K}$ is strong internal and $\varlimsup H_1\leq\varliminf H_2$ implies that ${\cal{K}}(H_1)\leq {\cal{K}}(H_1\cup H_2)\leq {\cal{K}}(H_2)$. 

\item A stronger form of monotonicity is \textbf{disjoint-monotone} that is if $H_1\cap H_2=\emptyset,{\cal{K}}(H_1)\leq{\cal{K}}(H_2)$ then ${\cal{K}}(H_1)\leq{\cal{K}}(H_1\cup H_2)\leq {\cal{K}}(H_2)$.

\item ${\cal{K}}$ is \textbf{part-shift-monotone} if $x>0,H_1\cap H_2=H_1\cap H_2+x=\emptyset$ then ${\cal{K}}(H_1\cup H_2)\leq {\cal{K}}(H_1\cup H_2+x)$.

\item ${\cal{K}}$ is \textbf{mean-monotone} if $H,K_1,K_2\in Dom({\cal{K}}),\sup K_1\leq{\cal{K}}(H)\leq\inf K_2$ implies that ${\cal{K}}(H\cup K_1)\leq{\cal{K}}(H)\leq{\cal{K}}(H\cup K_2)$. 

${\cal{K}}$ is \textbf{strong mean-monotone} if $H,K_1,K_2\in Dom({\cal{K}}),\varlimsup K_1\leq{\cal{K}}(H)\leq\varliminf K_2$ implies that ${\cal{K}}(H\cup K_1)\leq{\cal{K}}(H)\leq{\cal{K}}(H\cup K_2)$.

\item ${\cal{K}}$ is \textbf{base-monotone} if $H_1,H_2\in Dom({\cal{K}}), H_1\cap H_2=\emptyset$ then $$\min\{{\cal{K}}(H_1),{\cal{K}}(H_2)\}\leq{\cal{K}}(H_1\cup H_2)\leq\max\{{\cal{K}}(H_1),{\cal{K}}(H_2)\}.$$

${\cal{K}}$ is \textbf{countable-base-monotone} if $H_i\in Dom({\cal{K}})\ (i\in\mathbb{N}),\cup_1^{\infty}H_i\in Dom({\cal{K}}), H_i\cap H_j=\emptyset\ (\ i\ne j)$ then $$\inf\{{\cal{K}}(H_i):i\in\mathbb{N}\}\leq{\cal{K}}(\cup_{i\in\mathbb{N}}H_i)\leq\sup\{{\cal{K}}(H_i):i\in\mathbb{N}\}.$$

\item ${\cal{K}}$ is \textbf{union-monotone} if $B\cap C=\emptyset$

${\cal{K}}(A)\leq{\cal{K}}(A\cup B),{\cal{K}}(A)\leq{\cal{K}}(A\cup C)$ implies ${\cal{K}}(A)\leq{\cal{K}}(A\cup B\cup C)$ 

and 

${\cal{K}}(A\cup B)\leq{\cal{K}}(A),{\cal{K}}(A\cup C)\leq {\cal{K}}(A)$ implies ${\cal{K}}(A\cup B\cup C)\leq{\cal{K}}(A)$.

Moreover if any of the inequalities on the left hand side is strict then so is the inequality on the right hand side.

\item ${\cal{K}}$ is \textbf{d-monotone} if $L,B\in Dom({\cal{K}}),L\cap B=(L\cup B)\cap (B+x)=\emptyset$ then

${\cal{K}}(L)<{\cal{K}}(L\cup B), x>0$ implies ${\cal{K}}(L\cup B)<{\cal{K}}(L\cup B\cup B+x)$

and

${\cal{K}}(L)>{\cal{K}}(L\cup B), x<0$ implies ${\cal{K}}(L\cup B)>{\cal{K}}(L\cup B\cup B+x)$
\end{itemize}

\begin{prp}If ${\cal{K}}$ is (strong) internal, (strong) mean monotone then (strong) monotone.
\end{prp}
\P Let $\sup H_1\leq\inf H_2$. Set $K_1=\inf H_1, H=H_1, K_2=H_2$. Then we get that ${\cal{K}}(H_1)\leq{\cal{K}}(H_1\cup H_2)$. The other inequality and the "strong" part can be shown similarly.
\Pes

\begin{prp}If ${\cal{K}}$ is (strong) internal, finite-(countable)-independent, base-monotone then (strong) monotone. 
\end{prp}
\P If $\sup H_1\leq\inf H_2$ then by finite-independence we can assume that $H_1\cap H_2=\emptyset$. Then ${\cal{K}}(H_1)=\min\{{\cal{K}}(H_1),{\cal{K}}(H_2)\}\leq{\cal{K}}(H_1\cup H_2)\leq\max\{{\cal{K}}(H_1),{\cal{K}}(H_2)\}={\cal{K}}(H_2)$ by internality.

The "strong" part can be shown similarly.
\Pes

\begin{prp}If ${\cal{K}}$ is internal, finite-independent, disjoint-monotone then monotone. If ${\cal{K}}$ is internal, countable-independent, disjoint-monotone then strong-monotone. 
\end{prp}
\P Let $l=\sup H_1\leq\inf H_2$. By finite-independence we can assume that $l\notin H_1,H_2$. Then $H_1,H_2$ are disjoint and ${\cal{K}}(H_1)\leq\l\leq{\cal{K}}(H_2)$. Hence disjoint-monotonicity gives monotonicity. The "strong" part is similar.
\Pes

\begin{prp}If ${\cal{K}}$ is (strong) internal, disjoint-monotone then (strong) mean-monotone. 
\end{prp}
\P Let $H,K_1,K_2\in Dom({\cal{K}}),\sup K_1\leq{\cal{K}}(H)\leq\inf K_2$. Then $\sup K_1-H\leq{\cal{K}}(H)\leq\inf K_2-H$ and we get that ${\cal{K}}(K_1-H)\leq{\cal{K}}(H)\leq{\cal{K}}(K_2-H)$. Therefore ${\cal{K}}(H\cup K_1)={\cal{K}}(H\cup (K_1-H))\leq{\cal{K}}(H)\leq{\cal{K}}(H\cup (K_2-H))={\cal{K}}(H\cup K_2)$.

The "strong" version is similar.
\Pes

\begin{ex}In the definition of part-shift-monotonicity we cannot weaken the disjointness conditions significantly.

Let $H_1=[0,1],H_2=[0,1]\cup[11,12],x=1$. Then $Avg(H_1\cup H_2)=6$ while $Avg(H_1\cup H_2+1)=\frac{29}{6}<5$.
\end{ex}

\begin{ex}In the definition of base-monotonicity we cannot weaken the disjointness significantly.

Let ${\cal{K}}=Avg,H_1=[1,2]\cup[3,4],H_2=[1-\epsilon,1]\cup[3,4]$. $\epsilon$ can be chosen such that $Avg(H_2)>2.5$. Then obviously $\min\{{\cal{K}}(H_1),{\cal{K}}(H_2)\}={\cal{K}}(H_1)>{\cal{K}}(H_1\cup H_2)$.

\end{ex}

\begin{ex}In the definition of union-monotonicity we cannot weaken the disjointness significantly.

Let $A=[5,13], B=[2,5]\cup[13,17], C=[0,2]\cup[13,17]$. Then $Avg(A)=9,Avg(A\cup B)=9.5,Avg(A\cup C)=\frac{268}{28}>9$ while $Avg(A\cup B\cup C)=8.5$.
\end{ex}

\begin{prp}$Avg$ is part-shift-monotone. 
\end{prp}
\P Let $x>0,H_1\cap H_2=H_1\cap H_2+x=\emptyset$. There is anything to prove when both $H_1,H_2$ are $s$-sets for the same $s$. Then $Avg(H_1\cup H_2)=\frac{\mu^s(H_1)Avg(H_1)+\mu^s(H_2)Avg(H_2)}{\mu^s(H_1)+\mu^s(H_2)}<\frac{\mu^s(H_1)Avg(H_1)+\mu^s(H_2)Avg(H_2+x)}{\mu^s(H_1)+\mu^s(H_2)}=Avg(H_1\cup H_2+x)$
\Pes

\begin{prp}$Avg$ is disjoint-monotone. 
\end{prp}
\P Let $H_1\cap H_2=\emptyset, Avg(H_1)\leq Avg(H_2), H_1$ be an $s_1$-set, $H_2$ be an $s_2$-set. We have 3 cases.

1. $s_1=s_2=s$. Then $Avg(H_1\cup H_2)=\frac{\mu^s(H_1)Avg(H_1)+\mu^s(H_2)Avg(H_2)}{\mu^s(H_1\cup H_2)}\leq\frac{\mu^s(H_1)Avg(H_2)+\mu^s(H_2)Avg(H_2)}{\mu^s(H_1)+\mu^s(H_2)}=Avg(H_2)$. The other inequality is similar.

2. $s_1<s_2$. Then $Avg(H_1\cup H_2)=Avg(H_2)$.

3. $s_2<s_1$. Then $Avg(H_1\cup H_2)=Avg(H_1)$.
\Pes

\begin{prp}${\cal{M}}^{iso},{\cal{M}}^{acc},{\cal{M}}^{lis}$ are disjoint-monotone. \Pes
\end{prp}

\begin{prp}If ${\cal{K}}$ is disjoint-monotone, $H_1\cap H_2=\emptyset,{\cal{K}}(H_1)={\cal{K}}(H_2)=k$ then ${\cal{K}}(H_1\cup H_2)=k=\frac{{\cal{K}}(H_1)+{\cal{K}}(H_2)}{2}$. \Pes 
\end{prp}

\begin{prp}$Avg$ is d-monotone. 
\end{prp}
\P Let $Avg(H)<Avg(H\cup B)$. This implies that $H,B$ are $s$-sets for the same $s$. We can assume that $H\cap B=\emptyset$. Let $(H\cup B)\cap(B+x)=\emptyset,x>0$.

Then $Avg(H)<\frac{\mu^s(H)}{\mu^s(H)+\mu^s(B)}Avg(H)+\frac{\mu^s(B)}{\mu^s(H)+\mu^s(B)}Avg(B)$ implies that $Avg(H)<Avg(B)$.

Let us calculate $Avg(H\cup B\cup B+x)$. It is $\frac{\mu^s(H)Avg(H)+\mu^s(B)Avg(B)+\mu^s(B)(Avg(B)+x)}{\mu^s(H)+\mu^s(B)+\mu^s(B)}
=\frac{\mu^s(H)Avg(H)+2\mu^s(B)Avg(B)+\mu^s(B)x}{\mu^s(H)+2\mu^s(B)}>\frac{\mu^s(H)Avg(H)+2\mu^s(B)Avg(H)+\mu^s(B)x}{\mu^s(H)+2\mu^s(B)}>Avg(H)$.

The opposite inequality is similar.
\Pes

\begin{prp}$Avg$ is union-monotone.
\end{prp}
\P Let $B\cap C=\emptyset,Avg(A)\leq Avg(A\cup B),Avg(A)\leq Avg(A\cup C)$. Obviously we can assume that $A\cap B=A\cap C=\emptyset$ and $A,B,C$ are $s$-sets for the same $s$. Let $\mu=\mu^s$. We know that $Avg(A)\leq\frac{\mu(A)Avg(A)+\mu(B)Avg(B)}{\mu(A)+\mu(B)}$ and $Avg(A)\leq\frac{\mu(A)Avg(A)+\mu(C)Avg(C)}{\mu(A)+\mu(C)}$. 

Then $Avg(A\cup B\cup C)=\frac{\mu(A)Avg(A)+\mu(B)Avg(B)+\mu(C)Avg(C)}{\mu(A)+\mu(B)+\mu(C)}
\geq\frac{(\mu(A)+\mu(B))Avg(A)+\mu(C)Avg(C)}{\mu(A)+\mu(B)+\mu(C)}
=\frac{\mu(A)Avg(A)+\mu(C)Avg(C)+\mu(B)Avg(A)}{\mu(A)+\mu(B)+\mu(C)}
\geq\frac{(\mu(A)+\mu(C))Avg(A)+\mu(B)Avg(A)}{\mu(A)+\mu(B)+\mu(C)}=Avg(A)$

The opposite inequality is similar.
\Pes

\begin{prp}${\cal{M}}^{iso},{\cal{M}}^{acc}$ are union-monotone. \Pes
\end{prp}

\begin{prp}${\cal{M}}^{iso},{\cal{M}}^{acc},{\cal{M}}^{lis}$ are base-monotone. $Avg$ is countable-base-monotone. \Pes
\end{prp}

\begin{prp}Let ${\cal{K}}$ be union-monotone. 

(1) If $B\cap C=\emptyset$ and ${\cal{K}}(A)={\cal{K}}(A\cup B),{\cal{K}}(A)={\cal{K}}(A\cup C)$ then ${\cal{K}}(A)={\cal{K}}(A\cup B\cup C)$. 

(2) If $B\cap C=\emptyset$ and ${\cal{K}}(A)={\cal{K}}(A-B),{\cal{K}}(A)={\cal{K}}(A-C)$ then ${\cal{K}}(A)={\cal{K}}(A- (B\cup C))$. \Pes
\end{prp}

\begin{ex}${\cal{M}}^{eds}$ is not base-monotone. 

Let $H_1=\{\frac{1}{n},5+\frac{1}{n}:n\in\mathbb{N}\},H_2=\{1+\frac{1}{n},5+\frac{1}{n}+\frac{1}{2^n}:n\in\mathbb{N}\}$. Evidently $H_1\cap H_2=\emptyset$. Then ${\cal{M}}^{eds}(H_1)=2.5,{\cal{M}}^{eds}(H_2)=3$ hence $\min=2.5$. Using exactly the same method than in \cite{lamis} Example 10 one can show that ${\cal{M}}^{eds}(H_1\cup H_2)=\frac{0+1+5}{3}=2\not\geq 2.5$.
\end{ex}

\begin{prp}If ${\cal{K}}$ is base-monotone then it is convex.
\end{prp}
\P Let $I=[a,b],{\cal{K}}(H)\in I, L\subset I$. Then $H\cap(L-H)=\emptyset,\ {\cal{K}}(L-H)\in[a,b]$ which implies that $a\leq\min\{{\cal{K}}(H),{\cal{K}}(L-H)\}\leq{\cal{K}}(H\cup L)\leq\max\{{\cal{K}}(H),{\cal{K}}(L-H)\}\leq b$.
\Pes

\subsection{Continuity concepts}
\begin{itemize}
\item ${\cal{K}}$ is \textbf{slice-continuous} if $H\in Dom({\cal{K}})$ then $H^{+\varliminf H},H^{-\varlimsup H}\in Dom({\cal{K}})$ and $f(x)={\cal{K}}(H^{-x})$ and $g(x)={\cal{K}}(H^{+x})$ are continuous where $Dom( f)=\{x:H^{-x}\in Dom({\cal{K}})\},\ Dom(g)=\{x:H^{+x}\in Dom({\cal{K}})\}$.

\item ${\cal{K}}$ is \textbf{bi-slice-continuous} if $H\in Dom({\cal{K}})$ then $H^{+\varliminf H},H^{-\varlimsup H}\in Dom({\cal{K}})$ and 
$f(x,y)={\cal{K}}(H^{-x}\cup H^{+y})$ is continuous where $Dom(f)=\{(x,y):H^{-x}\cup H^{+y}\in Dom({\cal{K}})\}$.

\item ${\cal{K}}$ is \textbf{part-slice-continuous} if $H_1,H_2\in Dom({\cal{K}})$ then $H_2^{+\varliminf H},H_2^{-\varlimsup H}\in Dom({\cal{K}})$ and $f(x)={\cal{K}}(H_1\cup H_2^{-x})$ and $g(x)={\cal{K}}(H_1\cup H_2^{+x})$ are continuous where $Dom(f)=\{x:H_1\cup H_2^{-x}\in Dom({\cal{K}})\}, Dom(g)=\{x:H_1\cup H_2^{+x}\in Dom({\cal{K}})\}$.

\item Let $H\in Dom({\cal{K}}), x\in\mathbb{R}$. We say that ${\cal{K}}$ is point-continuous at $x$ regarding $H$ if $$\lim_{\epsilon\to 0+0}{\cal{K}}(H-S(x,\epsilon))={\cal{K}}(H).$$ 
We call ${\cal{K}}$ \textbf{point-continuous} if this holds for all $H\in Dom({\cal{K}})$ and $x\in\mathbb{R}$.

\item ${\cal{K}}$ is \textbf{Cantor-continuous} if $H_i\in Dom({\cal{K}}), H_{i+1}\subset H_i,\ \cap_{n=1}^{\infty}H_i\in Dom({\cal{K}})$ implies that ${\cal{K}}(H_i)\to {\cal{K}}(\cap_{n=1}^{\infty}H_i)$. 

\item Let ${\cal{F}}\subset\{f:\mathbb{R}\to\mathbb{R}\text{ continuous}\}$ equipped with the $\sup$ norm/topology. ${\cal{K}}$ is \textbf{f-continuous} with respect to ${\cal{F}}$ if $H\in Dom({\cal{K}}),\ f\in{\cal{F}}$ then $f(H)\in Dom({\cal{K}})$ and the function $F_H:{\cal{F}}\to\mathbb{R},\ F_H(f)={\cal{K}}(f(H))$ is continuous for all $H\in Dom({\cal{K}})$ i.e. $f_n\to f$ in sup norm implies that ${\cal{K}}(f_n(H))\to {\cal{K}}(f(H))$.

We get some special cases when ${\cal{F}}=$ translations / reflections / contractions / dilations.

\item ${\cal{K}}$ is $\lambda$\textbf{-continuous} if all Lebesgue measurable sets are in $Dom\ {\cal{K}}$ and moreover if $H\in Dom\ {\cal{K}},\ H$ is Lebesgue measurable, $0<\lambda(H)<+\infty$ then $\forall\epsilon>0$ there is $\delta>0$ such that $K\in Dom\ {\cal{K}},\ K$ Lebesgue measurable, $\lambda((H-K)\cup(K-H))<\delta$ implies that $$|{\cal{K}}(H)-{\cal{K}}(K)|<\epsilon.$$
\end{itemize}

Slice-continuity seems to be a very natural and somehow expectable concept however most of our previously defined means does not satisfy it as the next examples show.

\begin{ex}Let $H=\{\frac{1}{n},1+\frac{1}{n},2+\frac{1}{n}:n\in\mathbb{N}\}$. Then ${\cal{M}}^{acc}, {\cal{M}}^{iso}, {\cal{M}}^{lis}$ are not slice-continuous at $1$ and not point-continuous at $0,1,2$ (everywhere else they are point-continuous).

Let $H=\{\frac{1}{n},2+\frac{1}{2^n}:n\in\mathbb{N}\}$. According to \cite{lamis} we know that $LAvg(H)={\cal{M}}^{eds}(H)=0$ which gives that $LAvg,{\cal{M}}^{eds}$ are not slice-continuous and not point-continuous at $0$ (everywhere else they are point-continuous).
\end{ex}

\begin{ex}$Avg$ is not slice-continuous.

Let $H=C\cup[1,2]$ where $C$ is the Cantor-set. Let $x=1$. Then $Avg(H^{-x})=Avg(C)=Avg^s(C)=\frac{1}{2}$ where $s=\frac{\log 2}{\log 3}$. But $\lim_{\epsilon\to 0+0}Avg(H^{-(x+\epsilon)})=\lim_{\epsilon\to 0+0}Avg^1(H^{-(x+\epsilon)})=\lim_{\epsilon\to 0+0}1+\frac{\epsilon}{2}=1$.
\end{ex}

\begin{prp}\label{pavgsbsc}$Avg^s\ (0< s\leq 1)$ is bi-slice-continuous.
\end{prp}
\P Set $\mu=\mu^s$.

Let $y\in\mathbb{R}$ such that $H^{-y}\in Dom(Avg^s)$. For given $\epsilon>0$ one can find $\delta>0$ such that $|z-y|<\delta$ and $H^{-z}\in Dom(Avg^s)$ implies that 
$$\left\lvert\frac{\int\limits_{H^{-z}}x d\mu}{\mu(H^{-z})}-\frac{\int\limits_{H^{-y}}x d\mu}{\mu(H^{-y})}\right\lvert=
\left\lvert\frac{\mu(H^{-y})\Big(\int\limits_{H^{-z}}x d\mu-\int\limits_{H^{-y}}x d\mu\Big) + (\mu(H^{-y})-\mu(H^{-z}))\int\limits_{H^{-y}}x d\mu}{\mu(H^{-z})\mu(H^{-y})}\right\lvert\leq$$
$$\left\lvert\frac{\int\limits_{H^{-z}-H^{-y}}x d\mu}{\mu(H^{-z})}\right\vert + \left\vert\frac{\mu(H^{-z}-H^{-y})\int\limits_{H^{-y}}x d\mu}{\mu(H^{-z})\mu(H^{-y})}\right\lvert<\epsilon$$
using that $z\mapsto\mu(H^{-z})$ is continuous.

The other case ($+y$) can be handled similarly. And so can both cases together which gives  bi-slice-continuity.
\Pes

\begin{rem}Obviously it does not hold for $s=0$ i.e. for $Avg^0=\A$.
\end{rem}

\begin{prp}\label{bscimppc}If ${\cal{K}}$  is bi-slice-continuous then it is point-continuous as well.
\end{prp}
\P By definition $f(x,y)={\cal{K}}(H^{-x}\cup H^{+y})$ is continuous at $(x,x)$.
\Pe

Obviously for finite sets $Avg$ is not point-continuous but apart from those it is as the next proposition shows.

\begin{prp}$Avg$ is point-continuous assuming that $H\in Dom(Avg)$ is not finite.
\end{prp}
\P Assume that $H\in Dom(Avg)$ is infinite. Let $x\in\mathbb{R}$. Then set $s=\lim_{\epsilon\to 0+0}\dim H-S(x,\epsilon)$ where $\dim$ is the Hausdorff dimension. The limit exists because $\epsilon_1<\epsilon_2$ implies that $\dim H-S(x,\epsilon_1)>\dim H-S(x,\epsilon_2)$. Note that $s>0$. Let us observe that $\dim H=s$ because of countable stability of the dimension function. Now there are 2 cases:

1. $\exists\epsilon_0>0$ such that $\epsilon>\epsilon_0$ implies that $\dim H-S(x,\epsilon)=s$. Then the statement follows from \ref{pavgsbsc} and \ref{bscimppc}.

2. No such $\epsilon_0>0$ exists. We show that this case cannot happen. Clearly $H-{x}=\bigcup_{n=1}^{\infty}H-S(x,\frac{1}{n})$. Because $\dim H-S(x,\frac{1}{n})<s$ we get that $\mu^s(H-S(x,\frac{1}{n}))=0$ which gives that $\mu^s(H)=0$ which would yield that $H\not\in Dom(Avg)$.
\Pes

\begin{prp}\label{pfiscmc}Let ${\cal{K}}$ be finite-independent and slice-continuous. Then ${\cal{K}}(H)={\cal{K}}(H\cap[\varliminf H,\varlimsup H])$.
\end{prp}
\P Let us observe that if $\varliminf H=\varlimsup H$ then the statement obviously holds. Now suppose that $\varliminf H\ne\varlimsup H$.

We know that $g(x)={\cal{K}}(H^{+x})$ is continuous. If $x<\varliminf H$ then by finite-independence $g(x)={\cal{K}}(H)$. Hence $g(\varliminf H)={\cal{K}}(H)$ but $g(\varliminf H)={\cal{K}}(H\cap[\varliminf H,+\infty))$. 

Let $H_1=H\cap[\varliminf H,+\infty)$. Applying similar argument for $\varlimsup H_1=\varlimsup H$ we get that ${\cal{K}}(H_1)={\cal{K}}(H_1\cap(-\infty,\varlimsup H_1])$.
\Pes

\begin{prp}Let ${\cal{K}}$ be finite-independent, slice-continuous and $d\in(\varliminf H,\varlimsup H)$. Then there are $x,y\in[\varliminf H,\varlimsup H]$ such that ${\cal{K}}(H\cap[x,y])=d$.
\end{prp}
\P If $d={\cal{K}}(H)$ then by \ref{pfiscmc} we are done ($x=\varliminf H,y=\varlimsup H$).

Suppose $d>{\cal{K}}(H)$ (the other case can be handled similarly). Then $f(x)={\cal{K}}(H\cap[x,\varlimsup H])$ is continuous, $f(\varliminf H)<d, f(\frac{d+\varlimsup H}{2})>d$ hence there is an $y$ with the required property.
\Pes

\begin{prp}Let ${\cal{K}}$ be finite-independent and slice-continuous. If $H\in Dom({\cal{K}}),\ \varliminf H<{\cal{K}}(H)<\varlimsup H$ then there is $x\in(\varliminf H,\varlimsup H)$ such that $$\frac{{\cal{K}}(H^{-x})+{\cal{K}}(H^{+x})}{2}={\cal{K}}(H).$$
\end{prp}
\P Observe that $\lim_{x\to\varliminf H+0}{\cal{K}}(H^{-x})=\varliminf H, \lim_{x\to\varliminf H+0}{\cal{K}}(H^{+x})={\cal{K}}(H)$, and $\lim_{x\to\varlimsup H-0}{\cal{K}}(H^{+x})=\varlimsup H, \lim_{x\to\varlimsup H-0}{\cal{K}}(H^{-x})={\cal{K}}(H).$ That gives that
$$\lim_{x\to\varliminf H+0}\frac{{\cal{K}}(H^{-x})+{\cal{K}}(H^{+x})}{2}<{\cal{K}}(H),
\lim_{x\to\varlimsup H-0}\frac{{\cal{K}}(H^{-x})+{\cal{K}}(H^{+x})}{2}>{\cal{K}}(H).$$
Then refering to the continuity of $\frac{{\cal{K}}(H^{-x})+{\cal{K}}(H^{+x})}{2}$ completes the proof.
\Pes

\begin{prp}If ${\cal{K}}$ is slice-continuous then there is $x\in[\varliminf H,\varlimsup H]$ such that $\frac{{\cal{K}}(H^{-x})+{\cal{K}}(H^{+x})}{2}=x$.
\end{prp}
\P Let $f(x)=\frac{{\cal{K}}(H^{-x})+{\cal{K}}(H^{+x})}{2}-x$. By assumption $f$ is continuous. There is $\epsilon>0$ such that if $x\in (\varliminf H,\varliminf H+\epsilon)$ then $f(x)\geq 0$ while when $x\in(\varlimsup H-\epsilon,\varlimsup H)$ then $f(x)\leq 0$.
\Pes

\begin{prp}Let ${\cal{K}}$ be monotone, finite-independent and part-slice-continuous. Then ${\cal{K}}$ is strongly monotone.
\end{prp}
\P  Let $\varlimsup H_1\leq \varliminf H_2$ (the opposite case can be handled similarly). If $\varlimsup H_1<\varliminf H_2$ then we can assume that $\sup H_1\leq\inf H_2$ and by monotonicity we get the statement.

Let $\varlimsup H_1=\varliminf H_2=t$. We know that $g(x)={\cal{K}}(H_1\cup H_2^{+x})$ is continuous. Then $\lim_{x\to t,x<t}g(x)=g(t)$ but if $x<t$ then $g(x)={\cal{K}}(H_1\cup H_2)$ by finite-independence hence ${\cal{K}}(H_1\cup H_2^{+t})={\cal{K}}(H_1\cup H_2)$. Similarly using $f(x)={\cal{K}}(H_1^{-x}\cup H_2^{+t})$ we get that ${\cal{K}}(H_1^{-t}\cup H_2^{+t})={\cal{K}}(H_1\cup H_2^{+t})={\cal{K}}(H_1\cup H_2)$.

By monotonicity ${\cal{K}}(H_1^{-t})\leq{\cal{K}}(H_1^{-t}\cup H_2^{+t})$ but ${\cal{K}}(H_1^{-t})={\cal{K}}(H_1)$.
The other inequality for $H_2$ is similar.
\Pes

\begin{prp}Let $H_n\in Dom\ {\cal{K}}\ (n\in\mathbb{N}),\ H_n\to x\in\mathbb{R}$ in the sense that $\forall\epsilon>0\ \exists N\in\mathbb{N}$ such that $n>N$ implies that $H_n\subset(x-\epsilon,x+\epsilon)$. Then ${\cal{K}}(H_n)\to x$. \Pes
\end{prp}

\begin{prp}${\cal{M}}^{acc}$ is Cantor-continuous.
\end{prp}
\P For $H_i$ let $l_i\in\mathbb{N}$ such that $H_i^{(l_i)}\ne\emptyset,H_i^{(l_i+1)}=\emptyset$. By $H_{i+1}\subset H_i$ we get $l_{i+1}\leq l_i$ which gives that there is $k_0\in\mathbb{N}$ such that $i,j>k_0$ implies that $l_i=l_j$. We know that $H_i^{(l_i)}$ is finite hence there is $k_1>k_0$ such that $H_i^{(l_i)}=H_j^{(l_j)}$ when $i,j>k_1$ therefore ${\cal{M}}^{acc}(H_i)={\cal{M}}^{acc}(H_j)$.
\Pes

\begin{ex}\label{exancc}$Avg$ is not Cantor-continuous.

Let $H_{2n}=[0,\frac{1}{n}]\cup[1,\frac{2}{n}],H_{2n+1}=[0,\frac{2}{n}]\cup[1,\frac{1}{n}]\ (n\in\mathbb{N})$. Then $\cap_1^{\infty}H_n=\{0,1\},Avg(\cap_1^{\infty}H_n)=\frac{1}{2}$ but $Avg(H_n)$ is not even convergent because $Avg(H_{2n})\to\frac{1}{3},Avg(H_{2n+1})\to\frac{2}{3}$. \Pes
\end{ex}

\begin{prp}Let $H_n,H$ be bounded $s$-sets, $H=\cap_1^{\infty}H_n$. Then $Avg(H_n)\to Avg(H)$.
\end{prp}
\P $|Avg(H_n)-Avg(H)|=\left\lvert\frac{\int_{H_n}x d\mu^s}{\mu^s(H_n)}-\frac{\int_{H}x d\mu^s}{\mu^s(H)}\right\lvert=$
$$\left\lvert\frac{\mu^s(H)(\int_{H_n}x d\mu^s-\int_{H}x d\mu^s)+(\mu^s(H)-\mu^s(H_n))\int_{H}x d\mu^s}{\mu^s(H_n)\mu^s(H)}\right\lvert=$$
$$\left\lvert\frac{\mu^s(H)(\int_{H_n-H}x d\mu^s)+(\mu^s(H)-\mu^s(H_n))\int_{H}x d\mu^s}{\mu^s(H_n)\mu^s(H)}\right\lvert\to 0$$
using that $\mu^s(H_n)\to\mu^s(H)$.
\Pes

\begin{cor}If $0\leq s\leq 1$ then $Avg^s$ is Cantor-continuous. \Pes
\end{cor}

\begin{ex}${\cal{M}}^{lis}$ is not Cantor-continuous.

Let $H_{n}=\{1,1-\frac{1}{k},4,4+\frac{1}{k}:k\in\mathbb{N}, k\geq n\}\cup\{2,2-\frac{1}{k}:k\in\mathbb{N}\}$ Then ${\cal{M}}^{lis}(\cap_1^{\infty}H_n)=2$ while ${\cal{M}}^{lis}(H_n)=2.5\ (\forall n)$.  \Pes
\end{ex}

\begin{prp}${\cal{M}}^{lis}$ is f-continuous with respect to ${\cal{F}}=\{f:\mathbb{R}\to\mathbb{R}\text{ continuous}\}$.
\end{prp}
\P Obviously ${\cal{M}}^{lis}(f(H))=\frac{\min f(H)'+\max f(H)'}{2}$.
Now the statement simply follows from the following three facts: 

0.  If $f\in{\cal{F}},c\in\mathbb{R},g(x)=f(x)+c, H$ is compact then $c+\min f(H)=\min g(H),c+\max f(H)=\max g(H)$. 

1. If $f,g\in{\cal{F}},f\leq g, H\subset\mathbb{R}$ compact then $\min f(H)\leq\min g(H),\max f(H)\leq\max g(H)$.

2. If $f\in{\cal{F}}, H$ is bounded then $f(H')=f(H)'$.
\Pes

\begin{ex}$Avg$ is not f-continuous with respect to ${\cal{F}}=\{f:\mathbb{R}\to\mathbb{R}\text{ continuous}\}$.

Let $H=[0,1]\cup[2,3]$,
$$f(x)=\begin{cases}
0 & \text{if }x<1 \\
x-1 & \text{if }x\in[1,2] \\
1 & \text{if }x>2
\end{cases},$$
$$
f_{2n}(x)=\begin{cases}
0 & \text{if }x\in[-\infty,0] \\
\frac{x}{n} & \text{if }x\in[0,1] \\
(1-\frac{3}{2n})x+\frac{5}{2n}-1 & \text{if }x\in[1,2] \\
\frac{1}{2n}x+1-\frac{3}{2n} & \text{if }x\in[2,3] \\
1 & \text{if }x\in[3,+\infty) \\
\end{cases},$$
$$f_{2n+1}(x)=\begin{cases}
0 & \text{if }x\in[-\infty,0] \\
\frac{x}{2n} & \text{if }x\in[0,1] \\
(1-\frac{3}{2n})x+\frac{4}{2n}-1 & \text{if }x\in[1,2] \\
\frac{1}{n}x+1-\frac{3}{n} & \text{if }x\in[2,3] \\
1 & \text{if }x\in[3,+\infty) \\
\end{cases}.
$$
Obviously $f_n\to f$ in the $\sup$ norm, $f(H)=\{0,1\},f_{2n}(H)=[0,\frac{1}{n}]\cup[1-\frac{1}{2n},1],f_{2n+1}(H)=[0,\frac{1}{2n}]\cup[1-\frac{1}{n},1]$. Hence $Avg(f(H))=\frac{1}{2}$ but $Avg(f_n(H))$ is not even convergent similarly to Example \ref{exancc}.
\Pes
\end{ex}

One might think that $Avg^1$ might fulfill this property but the following example shows the contrary.

\begin{ex}$Avg^1$ is not f-continuous with respect to ${\cal{F}}=\{f:\mathbb{R}\to\mathbb{R}\text{ continuous}\}$.

Let $C$ denote the Cantor set, $f$ denote the Cantor function $(f:[0,1]\to[0,1])$. Let $(f_n)$ be the usual sequence of functions that approximate $f$ $(f_n:[0,1]\to[0,1])$. It is known that $f_n\to f$ in the $\sup$ norm.

Let $H=[0,\frac{1}{3}]\cup C$. For all $n\in\mathbb{N}\ f_n(H)=[0,\frac{1}{2}]\cup H_n$ where $H_n$ is of measure 0. Therefore $Avg^1(f_n(H))=\frac{1}{4}$. While $f(H)=f(C)=[0,1]$ which implies that $Avg^1(f(H))=\frac{1}{2}$.
\Pes
\end{ex}

\subsection{Other various properties}
\begin{itemize}
\item If ${\cal{K}}(H^*)={\cal{K}}(H)$ then ${\cal{K}}$ is called \textbf{condensed} where $H^*$ consists of the condensation points of $H$ where $x$ is a condensation point of $H$ if $\forall \epsilon>0\ |S(x,\epsilon)\cap H|>\aleph_0$. 

\item ${\cal{K}}$ is ${\cal{I}}$\textbf{-independent} for an ideal ${\cal{I}}$ if $H\in Dom({\cal{K}})-{\cal{I}}$ implies that $H\cup V\in Dom({\cal{K}}),\ {\cal{K}}(H)={\cal{K}}(H\cup V)$ where $V\in{\cal{I}}$ arbitrary. Recall that ${\cal{I}}\subset P(\mathbb{R})$ is an ideal if it is closed for finite union and it is descending. We emphasize that ${\cal{I}}$ is not necessarily a subset of $Dom({\cal{K}})$.

When ${\cal{I}}$ is the set of finite sets then we get back the notion of \textbf{finite-independent}.

Let us highlight the case when ${\cal{I}}$ is the set of countable sets when we call ${\cal{K}}$ \textbf{countable-independent}.

\item Let $H\in Dom({\cal{K}}), x\in\mathbb{R}$ such that $\varlimsup H\leq\varliminf H+x$ or $\varlimsup H+x\leq\varliminf H$. Then ${\cal{K}}$ is called \textbf{self-shift-invariant} if ${\cal{K}}(H\cup (H+x))={\cal{K}}(H)+\frac{x}{2}$.

\item  If $x\in\mathbb{R},H_1\cup H_2\in Dom({\cal{K}}),H_1\cup(H_2+x)\in Dom({\cal{K}}), H_1\cap H_2=H_1\cap H_2+x=\emptyset$ implies that $sign({\cal{K}}(H_1\cup(H_2+x))-{\cal{K}}(H_1\cup H_2))=sign(x)$ and $$|{\cal{K}}(H_1\cup(H_2+x))-{\cal{K}}(H_1\cup H_2)|\leq |x|$$ then ${\cal{K}}$ is called \textbf{part-shift-invariant}.

\end{itemize}

\begin{prp}$\cal{K}$ is finite independent iff ${\cal{K}}(H-V)={\cal{K}}(H)$ for any finite set $V$.
\end{prp}
\P $(\Rightarrow)$: $(H-V)\cup(H\cap V)=H$ and $H\cap V$ is finite.

$(\Leftarrow)$: $(H\cup V)-(V-H)=H$ and $V-H$ is finite.
\Pe

A similar argument shows:
\begin{prp}$\cal{K}$ is ${\cal{I}}$-independent (${\cal{I}}$ is an ideal) iff ${\cal{K}}(H-V)={\cal{K}}(H)$ for any set $V\in{\cal{I}}$.  \Pes
\end{prp}

\begin{prp}$\cal{K}$ is countable-independent then $\cal{K}$ is closed iff it is accumulated.  
\end{prp}
\P $cl(H)-H'$ consists of the isolated points of $H$ which is countable.
\Pe

$LAvg$ is closed but not accumulated hence it is not countable-independent.

\begin{ex}In the definition of part-shift-invariance we cannot weaken the disjointness significantly.

Let $H_1=[0,2]\cup[3,4], H_2=[2,3], x=1$. Then $Avg(H_1\cup H_2)=2$ while $Avg(H_1\cup H_2+x)=\frac{11}{6}$.
\end{ex}

\begin{prp}$Avg$ is condensed. \Pes
\end{prp}

\begin{prp}$Avg$ is part-shift-invariant.
\end{prp}
\P We can assume that $H_1,H_2$ are $s$-sets for the same $s$. Let $\mu=\mu^s$. $Avg(H_1\cup H_2+x)=\frac{\mu(H_1)Avg(H_1)+\mu(H_2)Avg(H_2+x)}{\mu(H_1)+\mu(H_2)}=\frac{\mu(H_1)Avg(H_1)+\mu(H_2)Avg(H_2)}{\mu(H_1)+\mu(H_2)}+\frac{x\mu(H_2)}{\mu(H_1)+\mu(H_2)}=Avg(H_1\cup H_2)+\frac{x\mu(H_2)}{\mu(H_1)+\mu(H_2)}$.
\Pes

\begin{prp}If $\cal{K}$ is countably independent then it is accumulated as well. Moreover it is condensed.
\end{prp}
\P The isolated points compose a countable set.

To prove the other statement observe that $|H-H^*|\leq\aleph_0$.
\Pes

\begin{prp}If ${\cal{K}}$ is accumulated and Cantor-continuous then it is condensed.
\end{prp}
\P Being accumulated implies that ${\cal{K}}(H)={\cal{K}}(H')={\cal{K}}(H^{(2)})=\dots$. We know that $H^*=\cap_1^{\infty}H^{(n)}$. Now Cantor-continuity gives that ${\cal{K}}(H^*)=\lim_{n\to\infty}{\cal{K}}(H^{(n)})={\cal{K}}(H)$.
\Pes


{\footnotesize

\smallskip
\noindent
Dennis G\'abor College, Hungary 1119 Budapest Fej\'er Lip\'ot u. 70.

\noindent 
E-mail: losonczi@gdf.hu, alosonczi1@gmail.com\\
}


\begin{thebibliography}{www}

\bibitem{bb} J. M. Borwein, P. B. Borwein, {\em The way of all means}, Amer. Math. Monthly {\em 94} (1987), 519–-522.

\bibitem{bullen} P. S. Bullen, {\em Handbook of means and their inequalities}, vol. 260 Kluwer Academic Publisher, Dordrecht, The Netherlands (2003).

\bibitem{darPal} Z. Dar\'oczy and Zs. P\'ales, {\em On functional equations involving means}, Publ. Math. Debrecen {\em 62} no. 3--4 (2003),
363--377.

\bibitem{be} B. Ebanks, {\em Looking for a few good means}, Amer. Math. Monthly {\em 119} (2012), 658--669.

\bibitem{hajja} M. Hajja, {\em Some elementary aspects of means},  International Journal of Mathematics and Mathematical Sciences, Means and Their Inequalities, Volume 2013, Article ID 698906, 1--9.

\bibitem{laesv} A. Losonczi, {\em Extending means to several variables}, arXiv prepint

\bibitem{lamis} A. Losonczi, {\em Means of infinite sets I}, arXiv prepint

\bibitem{lambm} A. Losonczi, {\em Measures by means, means by measures}, arXiv preprint


\end{thebibliography}
\end{document}